\documentclass[12pt]{article}

\usepackage{amsmath,amsthm,amssymb,latexsym, amsfonts}

\newtheorem{theorem}{Theorem}[section] 
\newtheorem{lemma}[theorem]{Lemma}     
\newtheorem{corollary}[theorem]{Corollary}

\title{Rational Distances with Rational Angles}
\author{Ryan Schwartz \and J\'ozsef Solymosi\footnote{The second author was supported by a Sloan Fellowship and NSERC and OTKA  grants.} \and Frank de Zeeuw \\ \\ Department of Mathematics \\ University of British Columbia \\ Vancouver, B.C., Canada V6T1Z2 \\ Email: \{ryano,solymosi,fdezeeuw\}@math.ubc.ca}

\begin{document}
\maketitle

\begin{abstract}
In 1946 Erd\H os asked for the maximum number of unit distances,
$u(n)$, among $n$ points in the plane.  He showed that $u(n)>
n^{1+c/\log\log n}$ and conjectured that this was the true magnitude.
The best known upper bound is $u(n)<cn^{4/3}$,  due to
Spencer, Szemer\'edi and Trotter.  We show that the upper bound
$n^{1+6/\sqrt{\log n}}$ holds if we only consider unit distances with rational
angle, by which we mean that the line through the pair of points makes
a rational angle in degrees with the $x$-axis.  Using an algebraic
theorem of Mann we get a uniform bound on the number of paths between
two fixed vertices in the unit distance graph, giving a contradiction
if there are too many unit distances with rational angle.  This bound
holds if we consider rational distances instead of unit distances as
long as there are no three points on a line.  A superlinear lower
bound is given, due to Erd\H os and Purdy.  If we have at most
$n^{\alpha}$ points on a line then we get the bound $O(n^{1+\alpha})$ or
$n^{1+\alpha+6/\sqrt{\log n}}$ for the number of rational distances with 
rational angle depending on whether $\alpha\ge 1/2$ or $\alpha < 1/2$ 
respectively.
\end{abstract}

\section{Introduction}
A famous problem of Erd\H os from 1946 \cite{Erdo46} concerns the
maximal number of unit distances among $n$ points in the plane; we
will denote this number by $u(n)$.  He showed that $u(n)>
n^{1+c/\log\log n}$, using a $\sqrt{n}\times \sqrt{n}$ piece of a
scaled integer lattice, and conjectured that this was the true
magnitude.  The best known upper bound is $u(n)<cn^{4/3}$, first
proved by Spencer, Szemer\'edi and Trotter in 1984 \cite{Spen84}.
This bound has several other proofs, the simplest of which was the
proof by Sz\'ekely \cite{Szek97}, using a lower bound for the crossing
number of graphs.  A recent result of Matou\v{s}ek \cite{Mato11} shows
that the number of unit distances is bounded above by $cn\log
n\log\log n$ for most norms.  As a general reference for work done on
the unit distances problem, see \cite{Bras06}.

We will show that the upper bound $n^{1+6/\sqrt{\log n}}$ holds if we only
consider unit distances that have \emph{rational angle}, by which we
mean that the line through the pair of points makes a rational angle
in degrees with the $x$-axis (or equivalently, its angle in radians,
divided by $\pi$, is rational).  Under this restriction, we can use an
algebraic theorem of Mann \cite{Mann65} to get a uniform bound on the
number of paths between two fixed vertices in the unit distance graph, which 
will lead to a contradiction if there are too many unit
distances with rational angle between the points.

In fact, our proof also shows that the bound $n^{1+6/\sqrt{\log n}}$ holds for
the number of rational distances with rational angles, if we
have no three points on a line.  The lower bound, $n^{1+c/\log\log n}$, of 
Erd\H{o}s does not apply in this case as we are restricted to rational angles.  
But a construction of Erd\H{o}s and Purdy gives a superlinear lower bound for 
unit (and hence rational) distances with rational angles.

If instead we allow up to $n^{\alpha}$ points on a line where $1/2\le\alpha\le 
1$, the number of rational distances with rational angles is bounded by $4 
n^{1+\alpha}$.  This bound is tight up to a constant factor with the lower bound 
now coming from an $n^{1-\alpha}\times n^{\alpha}$ square grid.  If we allow up 
to $n^{\alpha}$ points on a line where $0<\alpha<1/2$, the number of rational 
distances with rational angles is bounded above by $n^{1+\alpha+6/\sqrt{\log 
n}}$.  We get a lower bound of $cn^{1+\alpha}$ from $n^{1-\alpha}$ horizontal 
lines each containing $n^{\alpha}$ rational points so that no three points on 
different lines are collinear.

In Section~\ref{sec:main} we will state our main results and give an outline of 
the proof.  Section~\ref{sec:mann} contains the algebraic tools that we will 
use, including, for completeness, a proof of Mann's Theorem.  In 
Section~\ref{sec:graph} we use the bounds obtained from Mann's theorem and some 
graph theory to prove our main results.  In Section~\ref{sec:lower} we give 
lower bounds for the main results.

\section{Main Results and Proof Sketch} \label{sec:main}
We will say that a pair of points in $\mathbb{R}^2$ has \emph{rational angle}
if the line segment between them, viewed as a complex number $z = re^{\pi
i\gamma}$, has $\gamma\in\mathbb{Q}$.  Our first result is the following.  
\begin{theorem} \label{thm:unit}
Given $n$ points in $\mathbb{R}^2$, the number of pairs of points with unit
distance and rational angle is at most $n^{1+6/\sqrt{\log n}}$.
\end{theorem}

Roughly speaking, our proof goes as follows. Given $n$ points in the plane, we
construct a graph with the points as vertices, and as edges the unit line
segments that have rational angle.  We can represent these unit line segments as
complex numbers, which must be roots of unity because of the rational angle condition. 
Then if this graph has many edges, it should have many cycles of a given length
$k$, and each such cycle would give a solution to the equation $$\sum_{i=1}^k \zeta_i =0, $$
with $\zeta_i$ a root of unity. Using an algebraic theorem of Mann from 1965 \cite{Mann65}, 
we could give a uniform bound on the number of such solutions, depending only on
$k$ (under the non-degeneracy condition that no subsum vanishes).  If the number
of non-degenerate cycles goes to infinity with $n$, this would give a contradiction.

However, dealing with cycles of arbitrary length is not so easy, so instead in
our proof we count non-degenerate paths (which we will call \emph{irredundant} 
paths) of length $k$ between two fixed vertices, which correspond to solutions of
the equation $$\sum_{i=1}^k \zeta_i =a, $$ where $a\in \mathbb{C}, a\ne 0$ 
corresponds to the line segment between the two points. We have extended Mann's
theorem to this type of equation, giving a similar upper bound and proving our 
result.

In fact, in our proof it turns out that it is not necessary for the lengths to be
$1$, but that they only need to be rational. 
This is because our extension of Mann's theorem also works for equations of the type
$$\sum_{i=1}^k a_i \zeta_i =a, $$
where $a_i\in \mathbb{Q}$ and $a\in \mathbb{C}, a\ne 0$. This leads to the 
following results.

\begin{theorem} \label{thm:rat}
Suppose we have $n$ points in $\mathbb{R}^2$, no three of which are on
a line.  Then the number of pairs of points with rational distance
and rational angle is at most $n^{1+6/\sqrt{\log n}}$.
\end{theorem}
The constant $6$ in this theorem and the next is not optimal, but is the smallest integer
that followed directly from our proof.
 
\begin{theorem} \label{thm:rat3}
Suppose we have $n$ points in $\mathbb{R}^2$, with no more than $n^{\alpha}$ on
a line, where $0<\alpha<1/2$.  Then the number of pairs of points with rational
distance and rational angle is at most $n^{1+\alpha+6/\sqrt{\log n}}$.
\end{theorem}

\begin{theorem} \label{thm:rat2}
Suppose we have $n$ points in $\mathbb{R}^2$, with no more than $n^{\alpha}$ on
a line, where $1/2\le\alpha\le 1$.  Then the number of pairs of points with
rational distance and rational angle is at most $4n^{1+\alpha}$.
\end{theorem}

\section{Mann's Theorem} \label{sec:mann}

For completeness we provide a proof of Mann's Theorem.  We then 
prove the extension that we will need to prove the main
result in the next section.

\begin{theorem}[Mann]
Suppose we have 
$$\sum_{i=1}^ka_i\zeta_i =0,$$
with $a_i\in \mathbb{Q}$, the $\zeta_i$ roots of unity, and no subrelations
$\displaystyle\sum_{i\in I}a_i\zeta_i =0$ where $\emptyset\neq I \subsetneq [k]$. Then
$$\left(\zeta_i/\zeta_j\right)^m = 1$$
for all $i,j$, with $\displaystyle m = \mathop{\prod_{p\leq k}}_{p\ 
\mathrm{prime}} p$.
\end{theorem}
\begin{proof}
 We can assume that $\zeta_1=1$ and $a_1=1$, so that we have $1 +
\sum_{i=2}^ka_i\zeta_i =0$.  We take a minimal $m$ such that $\zeta_i^m=1$ for
each $i$.\\  We will show that $m$ must be squarefree, and that a prime $p$ that
divides $m$ must satisfy $p\leq k$. Together these prove the theorem.\\
 Let $p$ be a prime dividing $m$. Write $m = p^j\cdot m^*$ with $(p,m^*) = 1$,
and use that to factor each $\zeta_i$ as follows:
$$\zeta_i = \rho^{\sigma_i}\cdot \zeta_i^*,$$
with $\rho$ a primitive $p^j$th root of unity so $$\rho^{p^j}=1,~~~\left( 
\zeta_i^* \right)^{p^{j-1}m^*} = 1,~~~0\leq \sigma_i\leq p-1.$$
Now reorganize the equation as follows:
$$0 = 1+ \sum_{i=2}^ka_i\zeta_i = 1+ \sum_{l=0}^{p-1}\alpha_\ell\rho^\ell = 
f(\rho),$$ where the coefficients are of the form
$$\alpha_\ell = \sum_{i\in I_\ell} a_i\zeta_i^* \in 
\mathbb{Q}(\zeta_2^*,\dots,\zeta_k^*) = K,$$ with $I_\ell=\{i\in [k] : 
\sigma_i=\ell\}$.  So $f$ is a polynomial
over the field $K$ of degree $\leq p-1$ and $f(\rho)=0$. The polynomial $f$
is not identically zero, since that would give a subrelation containing
strictly fewer than $k$ terms. To see this, observe that we must have
$\sigma_i\geq 1$ for at least one $i$, otherwise $\zeta_i^{m/p} =1$ for each 
$i$,
contradicting the minimality of $m$.\\ But we can compute the degree
of $\rho$ over $K$ to be
$$\deg_K(\rho) = \frac{\phi(m)}{\phi(p^{j-1}m^*)} =
\frac{\phi(p^j)}{\phi(p^{j-1})} = \left\{\begin{array}{cl} p-1 &
\text{if }j=1 \\ p & \text{if } j>1. \end{array}\right.$$ This is a
contradiction unless $j=1$, which proves that $m$ is
squarefree.\\ Knowing that $m$ is squarefree, we have $m = p\cdot m^*$
with $(p,m^*) =1$, and
$$\zeta_i = \rho^{\sigma_i}\cdot \zeta_i^*,~~~\rho^p=1,~~~\left(
\zeta_i^* \right)^{m^*} = 1,~~~0\leq \sigma_i\leq p-1.$$ Still
$f(\rho)=0$ for $f(x)$ a polynomial over $K$, not identically
zero. But we know (\cite{Lang94}, Ch. VI.3) that the minimal
irreducible polynomial of $\rho$ over $K$ is $F(x) =
x^{p-1}+x^{p-2}+\dots +x+1$, hence we must have $f(x) = c F(x)$ for
some $c\in K$. In particular, $f$ has $p$ terms, which implies that our
original relation had at least $p$ terms, so $k\geq p$.
\end{proof}

\begin{theorem}
Suppose we have 
$$\sum_{i=1}^ka_i\zeta_i =a,~~~\sum_{j=1}^ka_j^*\zeta_j^* =a,$$ with
$a\in \mathbb{C}, a\ne 0$, $a_i\in \mathbb{Q}$, roots of unity $\zeta_i$, and
no subrelations $\displaystyle\sum_{i\in I}a_i\zeta_i =0$ or
$\displaystyle\sum_{j\in J}a_j^*\zeta_j^* =0$ where
$\emptyset\neq I \subsetneq [k]$ and $\emptyset\neq J
\subsetneq [k]$.  Then for any $\zeta_j^*$ there is a $\zeta_i$ such
that
$$\left(\zeta_j^*/\zeta_i\right)^m = 1$$
with $\displaystyle m = \mathop{\prod_{p\leq 2k}}_{p\ \mathrm{prime}} p$.
\end{theorem}
\begin{proof}
We have $\sum a_i\zeta_i = a = \sum a_j^*\zeta_j^*$, which gives the
single equation
\begin{equation} \label{eqn:sum1}
\sum_{i=1}^ka_i\zeta_i-\sum_{j=1}^ka_j^*\zeta_j^* =0.
\end{equation}
Mann's Theorem does not apply immediately, because there might be
subrelations. But we can break the equation up into minimal
subrelations
\begin{equation} \label{eqn:sum2}
\sum_{i\in I_\ell}a_i\zeta_i-\sum_{j\in I_\ell^*}a_j^*\zeta_j^* =0,
\end{equation}
where each $I_\ell\neq\emptyset$, $I_\ell^*\neq \emptyset$, and there are
no further subrelations.\\ Given $\zeta_j^*$, there is such a minimal
subrelation of length $\leq 2k$ in which it occurs, and which must
also contain some $\zeta_i$.  Applying Mann's Theorem to this equation
gives $\left(\zeta_j^*/\zeta_i\right)^m = 1$ with $\displaystyle m =
\mathop{\prod_{p\leq 2k}}_{p\ \mathrm{prime}} p$.
\end{proof}

Note that in the above proof we require $a\ne 0$.  If $a=0$ and there is no 
proper subrelation as in \eqref{eqn:sum2} then \eqref{eqn:sum1} still has the 
subrelations \[\sum_{i=1}^ka_i\zeta_i=0,~~~~~ \sum_{j=1}^ka_j^*\zeta_j^*=0,\] so we 
cannot use Mann's Theorem to get a relation between a $\zeta_i$ and $\zeta_j^*$.

For $a\in\mathbb{C}, a\ne 0, k\in\mathbb{Z}, k>0$ we define $Z_a^k$ to
be the set of $k$-tuples of roots of unity $(\zeta_1,\dots,\zeta_k)$ for which
there are $a_i\in\mathbb{Q}$ such that $\sum_{i=1}^ka_i\zeta_i=a$ with no
subrelations, i.e.: 
\[Z_a^k=\{(\zeta_1, \dots, \zeta_k)\ |\ \exists a_i \in
\mathbb{Q} : \sum_{i=1}^k a_i\zeta_i = a, \sum_{i\in I}a_i\zeta_i\ne 0\
\mathrm{for}\ \emptyset \ne I\subset [k]\}.\]

\begin{corollary} \label{cor:bound}
Let $\displaystyle C(k) = \mathop{\prod_{p\leq 2k}}_{p\ \mathrm{prime}}p$. Given 
$a\in \mathbb{C}, a\ne 0$, $|Z_a^k|\le (k\cdot C(k))^k$.
\end{corollary}
\begin{proof}
 Fix an element $(\zeta_1,\dots,\zeta_k)\in Z_a^k$ and let $m = C(k)$ and $M_i =
 \zeta_i^{-m}$ for $1\leq i \leq k$.  Then for $\zeta_j^*$ in any
 element of $Z_a^k$, we have an $i$ such that $M_i
 \left(\zeta_j^*\right)^m =1$.  In other words $\zeta_j^*$ is a
 solution of $M_i x^m = 1$.  Each of these $k$ equations has $m =
 C(k)$ solutions, hence there are at most $k\cdot m = k\cdot C(k)$ choices for each
 $\zeta_j^*$.
\end{proof}


\section{Rational Distances and Mann's Theorem} \label{sec:graph}

We are now in a position to prove the main results.
Suppose we have a graph $G=G(V,E)$ on $v(G)=n$ vertices and
$e(G)=cn^{1+\alpha}$ edges.  We will denote the minimum degree in $G$
by $\delta(G)$.  The following lemma assures us that we can remove low-degree 
vertices from our graph without greatly affecting the number of
edges. 
\begin{lemma}
\label{lem:minDeg}
Let $G$ be as above.  Then $G$ contains a subgraph $H$ with
$e(H)=(c/2)n^{1+\alpha}$ edges such that $\delta(H)\ge (c/2)n^{\alpha}$.
\end{lemma}
\begin{proof}
We iteratively remove vertices from $G$ of degree less than $(c/2)n^{\alpha}$.
Then, the resulting subgraph $H$ has $\delta(H)\ge (c/2)n^{\alpha}$ and we
removed fewer than $(c/2)n^{1+\alpha}$ edges so $H$ contains more than $(c/2)n^{1+\alpha}$ edges.
\end{proof}
Note that the subgraph $H$ constructed above contains at least 
$v(H)=(c/2)n^{\alpha}$ vertices.

Suppose we are given a path on $k$ edges $P_k=p_0p_1\dots p_k$.  We call this
path \emph{irredundant} if \[\sum_{i\in I}\overrightarrow{p_ip_{i+1}}\ne 0\] for
any $\emptyset\ne I \subset \{0,1,\dots,k-1\}$.

\begin{proof}[Proof of Theorem~\ref{thm:rat}]
Let $G$ be the graph with the $n$ points in the plane as vertices and
the rational distances with rational angles between pairs of points as
edges.  Suppose there are $n^{1+f(n)}$ such distances for some positive function 
$f$.  Then $e(G)=n^{1+f(n)}$.  We will count the number of
irredundant paths $P_k$ in $G$, for a fixed $k$ that we will choose later.  
By Lemma~\ref{lem:minDeg} we can
assume that $e(G)\ge (1/2)n^{1+f(n)}, v(G)\ge (1/2)n^{f(n)}$ and $\delta(G)\ge 
(1/2)n^{f(n)}$.

The number of irredundant paths $P_k$ starting at any vertex $v$ is at
least \[N=\prod_{\ell=0}^{k-1}(\delta(G)-2^{\ell}+1),\] since, if we have
constructed a subpath $P_\ell$ of $P_k$, then at most $2^{\ell}-1$ of the at 
least $\delta(G)$ continuations are forbidden.
Thus the total number of irredundant paths $P_k$ is at least 
\[\frac{nN}{2} \ge (n/2)\prod_{\ell=0}^{k-1}((1/2)n^{f(n)}-2^{\ell}+1) \ge 
\frac{n^{kf(n)+1}}{2^{2k+1}}\]
if $2^k \le (1/2)n^{f(n)}$, which is true as long as $k < f(n)\log n/\log 2$.  
It follows that there are two vertices $v$ and $w$ with at least
\[\frac{N}{n} \ge (1/n)\prod_{\ell=0}^{k-1}((1/2)n^{f(n)}-2^{\ell}+1) \ge\frac{n^{kf(n)-1}}{4^k}\]
irredundant paths $P_k$ between them.  We will call the set of 
these paths $\mathcal{P}_{vw}$, so that we have $|\mathcal{P}_{vw}|\ge n^{kf(n)-1}/4^k$.

Given $P_k\in\mathcal{P}_{vw}$, $P_k=p_0p_1\dots p_k$, consider the $k$-tuple
$(\zeta_1, \dots, \zeta_k)$ where $\zeta_i$ is the root of unity in the
direction from $p_{i-1}$ to $p_i$, i.e.
$\zeta_i=\overrightarrow{p_{i-1}p_i}/|\overrightarrow{p_{i-1}p_i}|$.  Note
that $(\zeta_1,\dots,\zeta_k)\in Z_a^k$, because $P_k$ is irredundant.  Since there
are no three points on a line, this process gives an injective map from
$\mathcal{P}_{vw}$ to $Z_a^k$ so $|\mathcal{P}_{vw}|\le (k\cdot C(k))^k$ by
Corollary \ref{cor:bound}.  Thus \[\frac{n^{kf(n)-1}}{4^k}\le (k\cdot C(k))^k 
\Longrightarrow n^{kf(n)-1} \le (4k\cdot C(k))^k.\]
But this gives \[e^{(kf(n)-1)\log n} \le e^{k\log (4k\cdot C(k))} 
\Longrightarrow f(n) \le \frac{\log (4k) + \log (C(k))}{\log n} + \frac{1}{k}.\]

The term $\log(C(k))$ is the log of the product of the primes less than or equal to 
$2k$.  This is a well known number-theoretic function called the Chebyshev 
function and denoted by $\vartheta$, specifically $\vartheta(2k)=\log(C(k))$.  
We use the following bound on $\vartheta$ (for a proof see \cite{Apos76}): 
\[\vartheta(x) < 4x\log 2 < 3x,\quad\mathrm{for}\ x\ge2.\]
This gives \[f(n) < \frac{\log(4k)+6k}{\log n} + \frac{1}{k} < \frac{7}{\log n}k 
+ \frac{1}{k}.\]  Let $k$ be an integer such that $f(n)\log n/18 < k < f(n)\log 
n /14$, (possible since otherwise $f(n) = O(1/\log n)$ giving $n^{f(n)} = 
O(1)$).  Then the condition that $k< f(n)\log n/\log 2$ is  clearly satisfied, 
and we get \[f(n) < \frac{7}{\log n}\cdot \frac{f(n)\log n}{14} + 
\frac{18}{f(n)\log n} \Longrightarrow f(n) < \frac{6}{\sqrt{\log n}}.\]  This 
completes the proof.
\end{proof}

\begin{proof}[Proof of Theorem~\ref{thm:unit}]
In the statement of Theorem~\ref{thm:unit}, the requirement that there
are no three points on a line is unnecessary.  This is because, from
any point, there is only one unit distance in any direction.  Thus
we can apply the same proof as in Theorem~\ref{thm:rat} to
Theorem~\ref{thm:unit} without having to worry about multiple points
on a line.  Thus we also have a proof of Theorem~\ref{thm:unit}.
\end{proof}

Consider a path $P_k=p_0p_1\dots p_k$.  If the distance from $p_{i-1}$ to
$p_i$ is less than the distance from $p_{i-1}$ to any vertex on the line
connecting $p_{i-1}$ and $p_i$ and not in $P_{i-1}=p_0p_1\dots p_{i-1}$ then
$P_k$ is called a \emph{shortest path}.

\begin{proof}[Proof of Theorem~\ref{thm:rat3}]
This proof is almost the same as the proof of Theorem~\ref{thm:rat} except that
instead of considering all irredundant paths $P_k$, we only consider shortest
irredundant paths.  Suppose there are $n^{1+\alpha+f(n)}$ edges in the rational
distance graph.  Since there are at most $n^{\alpha}$ points on a line, we
get that from any vertex $v$ there are at least
\[N=\prod_{\ell=0}^{k-1}\biggl(\frac{\delta(G)}{n^{\alpha}}-2^{\ell}+1\biggr)\ge
\frac{n^{k f(n)}}{4^k}\] shortest irredundant paths $P_k$, if $k < f(n)\log 
n/\log 2$.  For any two vertices $v,w$ let
$\mathcal{P}_{v,w}$ be the set of shortest irredundant paths $P_k$ between $v$
and $w$.  Then there are two vertices $v,w$ such
that the number of shortest irredundant paths between $v$ and $w$ is at least
\[|\mathcal{P}_{v,w}|\ge \frac{n^{kf(n)-1}}{4^k}.\]  By Mann's Theorem, since we 
are looking at shortest irredundant paths, $|\mathcal{P}_{v,w}|\le (k\cdot 
C(k))^k.$  Let $k$ be an integer such that $f(n)\log n/18 < k < f(n)\log n/14$.  
Then \[\frac{n^{kf(n)-1}}{4^k}\le (k\cdot C(k))^k \Longrightarrow f(n) < 
\frac{6}{\sqrt{\log n}}.\]
\end{proof}

\begin{proof}[Proof of Theorem~\ref{thm:rat2}]
Assume we have a configuration of $n$ points with at most $n^{\alpha}$ on a
line, $1/2\le\alpha\le1$, and $n^{1+\alpha + f(n)}$ rational distances
with rational angles, for some positive function $f(n)$.

The graph $G$ on these points has $e(G)=n^{1+\alpha+f(n)}$.  
By Lemma~\ref{lem:minDeg} we can assume that $e(G)\ge n^{1+\alpha+f(n)}/2$, 
$v(G)\ge{n^{\alpha+f(n)}/2}$ and $\delta(G)\ge n^{\alpha+f(n)}/2$. We now count 
irredundant paths $P_2$ of length 2.  Note that an irredundant path on two edges 
is just a noncollinear path.

For any vertex $v$, since we have at most $n^{\alpha}$
points on a line, $v$ is the midpoint of at
least \[N=\delta(G)(\delta(G)-n^{\alpha})\ge \frac{n^{2(\alpha+f(n))}}{8}\]
paths $P_2$ if $f(n) \ge \log 4/\log n$ (if $f(n) < \log 4/\log n$ then 
$n^{f(n)} < 4$, completing the proof.)  Thus there are two vertices $v$ and $w$ 
with at least $(1/8)n^{2(\alpha+f(n))-1}$ noncollinear paths $P_2$ between them.

But by Corollary~\ref{cor:bound} there is a constant number of directions from 
each of $v$ and $w$.  Since we are looking at noncollinear paths $P_2$, the
direction from $v$ and the direction from $w$ uniquely determine the
midpoint for a path $P_2$.  Thus there are at most $(k\cdot C(k))^k=144$ 
noncollinear paths $P_2$ between $v$ and $w$, since $k=2$.

Putting the upper and lower bounds together we get that
$n^{2(\alpha+f(n))-1}\le 2^73^2$.  This gives \[f(n)\le \frac{7\log 2+2\log 
3}{2\log n} + \frac{1}{2} - \alpha \le \frac{7\log 2 + 2\log 3}{2\log n} < 
\frac{4}{\log n},\]
since $\alpha \ge 1/2$.  But this gives $n^{f(n)} < 4$, completing the proof.
\end{proof}

\section{Lower Bounds} \label{sec:lower}

In this section we give lower bounds for the theorems given in 
Section~\ref{sec:main}.

The bounds in Theorems~\ref{thm:unit} and~\ref{thm:rat} are not far
from optimal as the following construction of Erd\H{o}s and Purdy
\cite{Erdo76} shows. 

Suppose we have $n$ points, no three on a line, with the maximum
possible number of unit distances with rational angles; we call this
number $f(n)$.  Consider these points as the set $\{z_1,\dots,z_n\}$
of complex numbers.  For any $a\in\mathbb{C}$ with $|a|=1$, $a\ne z_i-z_j$
for any $i \ne j$, the set $\{z_1,\dots,z_n,z_1+a,\dots,z_n+a\}$ contains at
least $2f(n)+n$ unit distances since there are $f(n)$ amongst each of the sets
$\{z_1,\dots,z_n\}$ and $\{z_1+a,\dots,z_n+a\}$ and $|z_i-(z_i+a)|=1$
for each $i$.  This new set may have three points on a line, but we
show that we can choose $a$ appropriately so this is not the case.

Consider a pair of points $z_i$ and $z_j$.  For each $z_k$, the set of
points $\{z_k+a : |a|=1\}$ intersects the line through $z_i$ and $z_j$
in at most two points.  So there are at most two values of $a$ that
will give three points on a line.  There are $\binom{n}{2}$ pairs of
points and $n$ choices for $z_k$ so there are at most
$2n\binom{n}{2}=n^2(n-1)$ values of $a$ that make a point $z_k+a$
collinear with two points $z_i$ and $z_j$.  Similarly we have
$n^2(n-1)$ values of $a$ that make a point $z_k$ collinear with two
points $z_i+a$ and $z_j+a$.  Thus there are only finitely many values
of $a$ that give three points on a line.  There are infinitely many
choices for $a$ so we are done.

This shows that $f(2n) \ge 2f(n) + n$ for $n > 2$ and clearly
$f(2)=1$.  From this we get that $f(2^k)\ge 2^{k-1}(k-1)=2^{k-1}\log_2
(2^{k-1})$.  Taking $2^k\le n < 2^{k+1}$ we get that $f(n)\ge cn\log
n$ for all $n$.  This construction gives a lower bound for
Theorems~\ref{thm:unit} and~\ref{thm:rat}.

The bound in Theorem~\ref{thm:rat3} is not far from optimal.  In fact we can get
a lower bound of $cn^{1+\alpha}$.  Consider $n^{1-\alpha}$ lines parallel to the 
$x$-axis, and choose $n^{\alpha}$ rational points on each line such that no 
three points on different lines are collinear (this can always be done since 
there are infinitely many rational points to choose from).  There are 
$cn^{2\alpha}$ rational distances on each horizontal line and $n^{1-\alpha}$ 
such lines giving at least $cn^{1+\alpha}$ rational distances with rational 
angles (all the angles are zero).

The bound in Theorem~\ref{thm:rat2} is tight up to a constant
factor as can be seen by considering an $n^{1-\alpha}\times n^{\alpha}$
square grid.  Then there are at least $cn^{2\alpha}$ rational distances on each 
of the $n^{1-\alpha}$ horizontal lines in the grid containing $n^{\alpha}$ 
points.  This gives at least $cn^{1+\alpha}$ rational distances with rational 
angles (the angles are all zero).

\section*{Acknowledgment}

The authors would like to thank Jirka Matou\v{s}ek for making us aware of
the lower bound for unit distances with no three points on a line.  We are 
indebted to the anonymous referee for helpful comments and suggestions.

\bibliographystyle{plain} \bibliography{references}

\begin{thebibliography}{1}

\bibitem{Apos76}
T.M. Apostol.
\newblock {\em Introduction to Analytic Number Theory}, chapter 4.5:
  Inequalities for $\pi(n)$ and $p_n$, pages 82--85.
\newblock Springer, 1976.

\bibitem{Bras06}
P.~Brass, W.~Moser, and J.~Pach.
\newblock {\em Research Problems in Discrete Geometry}, chapter 5: Distance
  Problems, pages 183--257.
\newblock Springer, 2006.

\bibitem{Erdo46}
P.~Erd{\H o}s.
\newblock On sets of distances of $n$ points.
\newblock {\em American Mathematical Monthly}, 53(5):248--250, May 1946.

\bibitem{Erdo76}
P.~Erd{\H o}s and G.B. Purdy.
\newblock Some extremal problems in geometry, {IV}.
\newblock In {\em Proceedings of The Seventh Southeastern Conference on
  Combinatorics, Graph Theory, and Computing}, pages 307--322, 1976.

\bibitem{Lang94}
S.~Lang.
\newblock {\em Algebra}.
\newblock Addison-Wesley, 3rd edition, 1994.

\bibitem{Mann65}
H.B. Mann.
\newblock On linear relations between roots of unity.
\newblock {\em Mathematika}, 12:107--117, 1965.

\bibitem{Mato11}
J.~Matou{\v s}ek.
\newblock The number of unit distances is almost linear for most norms.
\newblock {\em Advances in Mathematics}, 226:2618--2628, 2011.
\newblock http://arxiv.org/abs/1007.1095.

\bibitem{Spen84}
J.~Spencer, E.~Szemer{\' e}di, and W.~Trotter.
\newblock Unit distances in the {E}uclidean plane.
\newblock In B.~Bollobas, editor, {\em Graph Theory and Combinatorics:
  Proceedings of the Cambridge Combinatorial Conference, in Honour of {P}aul
  {E}rd{\H o}s}, pages 293--303. Academic Press, 1984.

\bibitem{Szek97}
L.~Sz{\' e}kely.
\newblock Crossing numbers and hard {E}rd{\H o}s problems in discrete geometry.
\newblock {\em Combinatorics, Probability and Computing}, 6(3):353--358,
  September 1997.

\end{thebibliography}
\end{document}